\documentclass[amsfont amscd ,12pt]{amsart}

\usepackage{amssymb}
\usepackage{amsmath}
\usepackage{amscd}

\begin{document}

\title{SPHERE  THEOREM FOR MANIFOLDS WITH POSITIVE
CURVATURE}

\author{  Bazanfar\'e Mahaman}
\thanks{2000 mathematics subject classification
53C20, 53C21 }

\date{}

        \begin{abstract}

In this paper, we  prove that, for any integer $n\ge
2,$
there exists an $\epsilon_{n} \ge 0$ so that
 if $M$ is an  n-dimensional complete manifold
 with sectional curvature
$ K_{M}\ge 1$ and if $M$ has conjugate
radius bigger than $\frac{\pi}{2} $ and contains
a geodesic loop of length \mbox{$2(\pi
-\epsilon_{n}),$} then $M$
is diffeomorphic to
the Euclidian unit sphere  $S^{n}.$
\end{abstract}
 \maketitle
\section{Introduction}
One of the fundamental problems in Riemannian geometry
is to determine 
the relation between the topology and the geometry of
a Riemannian
manifold. In this way the Toponogov's theorem and the
critical point 
theory play an important rule.
Let $M$ be a complete Riemannian manifold and fix a
point $p$ in $M$ 
and define $d_{p}(x) = d(p,x).$ A point $q\neq p$ is
called
 a critical point of $d_{p}$ or simply of the point
$p$ if, for any 
nonzero vector $v\in T_{q}M$,there exists a minimal
geodesic $\gamma$
 joining $q$ to $p$ such that the angle 
$(v,\gamma'(0)\le\frac{\pi}{2}.$ Suppose $M$ is an
n-dimensional complete Riemannian
 manifold with sectional curvature $K_{M}\ge 1.$ By
Myers' theorem the 
diameter of $M$ is bounded from above by $\pi.$
 In [Ch] Cheng showed that the maximal value $\pi$ is
attained if and 
only if $M$ is isometric to the standard sphere.
 It was proved by Grove and Shiohama [GS] that if
$K_{M}\ge 1$ and the 
diameter of $M$  $diam(M)> \frac{\pi}{2}$ then $M$ is
homeomorphic
 to a sphere.

 Hence the problem of removing homeomorphism to
diffeomorphism or 
finding conditions to guarantee the diffeomorphism is
 of particular interest.
 In [Xi3] C. Xia showed that if $K_{M}\ge 1 $ and the
conjugate radius 
of $M$ $\rho(M)>\pi/2$ and if $M$ contains a geodesic
loop of
 length $2\pi$ then $M$ is isometric to $S^{n}(1).$
 \subsection{Definition}

   \textit{Let M be an n-dimensional  Riemannian
manifold and $p$ be a  
point in M. Let Conj(p) denote the set of
   first conjugate  points to p on all geodesics
issuing from  p. The 
conjugate radius $\rho(p)$ of $M$  at $p$ is defined
as}
         \begin{equation*}\rho(p) =
   d(p,Conj(p))\qquad \textrm{ if } Conj(p) \ne
\emptyset
   \end{equation*} and $$\rho(p) = +\infty
\qquad\textrm{ if } Conj(p) 
= \emptyset$$
   \textit{Then the conjugate radius of $M$ is:}$$\rho
(M) = \inf_{x\in 
M}\rho (x).$$
Many interesting results have been proved by using the
critical points 
theory
 and Toponogov's theorem [C], [GS],  [Pe], [S], [Sh],
[SS], [Xi1],
 [Xi2], [Xi3].
 In otherewise J. Cheeger and T. Colding in [CC] have
proven the 
following
 \subsection*{Theorem A}

\textit{There exists a number, $\epsilon(n)>0$,
depending only on the 
integer $n$ such that, for any two
Riemannian manifolds $Z_{1}, Z_{2}$, if
$d_{GH}(Z_{1},Z_{2})< 
\epsilon(n),$ then $Z_{1}$ and $Z_{2}$ are
diffeomorphic.}
where $d_{GH}(Z_{1},Z_{2}$ denote the Gromov-Hausdroff
distance.

    The purpose of this paper is to prove the
following :
    \subsection{Theorem}\label{**}

   \textit{ For any $n\ge 2$
    there exists a positif constant $\epsilon(n)$
depending only  on $n$ 
such that for any $\epsilon \le \epsilon(n)$,
     if $M$ is an n-dimensional
    complete connected Riemannian manifold with
sectional
curvature $K_{M}\ge 1$ and conjugate radius
    $\rho(M)>\frac{\pi}{2}$ and if  $M$ contains
a geodesic loop of length
$2(\pi-\epsilon)$ then $M$
    is diffeomorphic to an n-dimensional unit sphere
$S^{n}(1).$}
\section{Proof}
Since $K_{M}\ge 1,$  $M$ is compact. Let $i(M)$ denote
the injectivity
radius of $M.$

By definition we have $$ i(M) = \inf_{x\in M}
d(x,C(x)).$$
Since  $M$ is compact and the function
$x\mapsto d\left(x, C(x)\right)$ is
continuous, there  exists $p\in M$ such that
$i(M) = d(p,C(p)).$ Since $C(p)$
is compact there exists  $q\in C(p)$ such that
$i(M) = d(p,q) = d\left(p,
C(p)\right).$ Then

a) either there exists a minimal geodesic
$\sigma$ joining $p$ to $q$
such that $q$ is a conjugate  point of  $p$
 or

 b) there exists two minimal geodesics
$\sigma_{1}$  and $\sigma_{2}$ joining
$p$ to $q$ such that
$\sigma'_{1}(l) = - \sigma'_{2}(l), \quad l =d(p,q).$
See [C].
\subsection{Lemma}

\textit{Let  M be an n-dimensional complete, connected
Riemannian
 manifold with sectional  curvature
$ K_{M}\ge 1$ and the  conjugate radius
$\rho (M) >\frac{\pi}{2},$ then $i(M)
>\frac{\pi}{2}.$}

\begin{center}Proof of the lemma\end{center}

  If  $a)$ holds, then  $i(M) = d(p,q) > \pi/2 .$

Suppose  $b)$ holds . Since
  $q\in C(p),$ we have $p\in C(q)$  and consequently
$d(p,q) = d(q,C(q)).$ This implies that
  $\sigma_{1}'(0) = -\sigma_{2}'(0).$
  Set $D(x) = \max_{y\in M}d(x,y).$ Then
  $$D(x)\ge \max_{y\in C(x)} d(x,y)\ge \rho(M) > \pi/2
.$$
  Since $M$ is  compact, there exist a point $y\in M$
such that
 $D(x) = d(x,y) > \pi/2 $
and  $y$ is  the unique farthest  point and the
critical one
 for the distance $d(x,.).$
  Set $A(x) = y$; thus we  define a continuous  map
$A: M\mapsto M $ \mbox{(see [Xi3].}
By the  Berger-Klingenberg  theorem,
    $M$ is homeomorphic to the unit  sphere $S^{n}(1)$
and since $A(x)\ne x$ for all
   $ x\in M$ the Brouwer fixed point  theorem sets
that the
degree of  $A$ is
 $(-1)^{n+1}$ and consequently     $A$ is surjective.
  Let  $r\in M$ the point so that  $p = A(r).$
 Hence \mbox{$d(p,r) > \pi/2.$}
   If $r = q$ then $i(M) = d(p,q) > \pi/2 $ otherwise
there
 exists a minimal  geodesic from $q$ to $r$,
   $\sigma_{3}$
    such that
$$\measuredangle(\sigma_{3}'(0),-\sigma_{1}'(l))
\le \pi/2\quad or \quad
    \measuredangle(\sigma_{3}'(0),-\sigma_{2}'(l))\le
\pi/2.$$

    Suppose
$\measuredangle(\sigma_{3}'(0),-\sigma_{1}'(l))\le
\pi/2.$
   Applying the Toponogov's  theorem [T] to the hinge
$(\sigma_{1},\sigma_{3}),$ we have:
    \begin{eqnarray*}\cos d(p,r)\ge \cos d(p,q)\cos
d(q,r) + \sin 
d(p,q)\sin d(q,r)\cos\angle(
    \sigma_{3}'(0),-\sigma_{1}'(l))\end{eqnarray*}
     \begin{eqnarray}\label{2.33}
    \ge cos d(p,q)\cos d(q,r).\end{eqnarray}
     Since  $r$ is far from $p$ in the sense that
  $d(p,r) >\pi/2 $ then $r$ is near to $q$ i.e
$d(q,r)< \pi/2$  and from
    (\ref{2.33}) we have
 $$\cos d(p,q) < 0$$ and  consequently $$i(M)= d(p,q)
>\pi/2 $$ which 
proves the lemma.
    \subsection{Lemma}

    \textit{Let  $M$ be a complete connected
$n$-dimensional Riemannian 
manifold
with sectional  curvature $ K_{M}\ge 1$ and conjugate
radius
 $\rho (M) >\frac{\pi}{2} $. If $M$ contains a
geodesic loop of length
at least
$2(\pi-\epsilon)$ then $diam(M)\ge \pi-\tau(\epsilon)$
where
 $\tau(\epsilon)\mapsto 0$
 when $\epsilon\mapsto 0.$}

 \begin{center}Proof\end{center}

 Since $i(M)>\pi/2$ then there exist  $\delta >0$ such
that
 $i(M) > \pi/2 +\delta.$
    Let $\gamma$ be a loop with length
 $2\pi-2\epsilon.$
    Let
 $x = \gamma (0) =\gamma (2\pi-2\epsilon)$,
$y = \gamma(\pi/2  +\delta) $, $m = \gamma(\pi-\epsilon)$
    and $z = \gamma (\frac{3(\pi-\epsilon)}{2}-\delta)$
    
 Let 
 $$\gamma_{1} = \gamma_{/[0, \frac{\pi}{2}
+\delta]},\quad
\gamma_{2} = \gamma_{/[\frac{\pi}{2}+\delta, \pi
-\epsilon]}\quad\gamma_{3} = \gamma_{/[\pi -\epsilon, \frac{3(\pi
-\epsilon)}
{2}-\delta]} $$ 
and $\gamma_{4} =
\gamma_{/[\frac{3(\pi -\epsilon)}
{2}-\delta, 2\pi-2\epsilon]}.$

Then the geodesics $\gamma_{i} $ are minimal. Let
$\sigma$ be a minimal
geodesic joining $m$ and $x.$

We claim  that $L(\sigma) \ge \pi-\tau(\epsilon).$

Set $\alpha = \measuredangle
(\sigma'(0),-\gamma'(\pi -\epsilon))$ and $\beta =
\measuredangle
(\sigma'(0),\gamma'(\pi -\epsilon)).$
Applying the Toponogov's theorem to the
 triangles $(\gamma_{1}, \gamma_{2},\sigma)$ and   
$(\gamma_{3},
\gamma_{4},\sigma)$ respectively,  one can take two 
triangles
$(\overline{\gamma}_{1},
\overline{\gamma}_{2},\overline{\sigma})$  and
$(\overline{\gamma}_{3},
\overline{\gamma}_{4},\overline{\sigma})$ on 
the unit sphere
$S^{2}(1)$ with vertices
$\overline{x},\overline{y},\overline{m}$ and
 $\overline{x},\overline{z},\overline{m}$ respectively
satisfying:
$$L(\overline{\gamma}_{i}) = L(\gamma_{i}), \quad i
=1,2,3,4;
L(\overline{\sigma}) = L(\sigma);$$  hence
 $\overline{\alpha} \le \alpha,
\overline{\beta}\le \beta$ where  $\overline{\alpha}$
and
 $\overline{\beta}$ are
the angles at $\overline{m} $ of the triangles
 $(\overline{\gamma}_{1},
\overline{\gamma}_{2},\overline{\sigma})$  and
$(\overline{\gamma}_{3},
\overline{\gamma}_{4},\overline{\sigma})$
respectively.
 We have:$\alpha\le \pi/2$ or
$\beta\le \pi/2.$ Suppose, without lost the
generality, that
$\alpha\le \pi/2.$ Let
 $\overline{x}'$ be the antipodal point of  
$\overline{x}$ on the 
sphere
 $S^{2}(1)$  and
$\overline{\sigma}_{1}$ the minimal geodesic from
$\overline{m}$ to 
$\overline{x}'.$ If   $\overline{\alpha}'$ and
 $\overline{\beta}'$ are the angles at
$\overline{m},$  of triangles 
$(\overline{y},\overline{m},\overline{x}')$ and $
(\overline{z},\overline{m},\overline{x}')$
respectively then we have:
\begin{eqnarray} d((\overline{y},\overline{x}') = \pi
-d(\overline{y},\overline{x})\nonumber\\
= \pi - d(x,y)
=\frac{\pi}{2}  - \delta
\end{eqnarray}

Hence, using the trigonometric law on the triangle
$(\overline{y},\overline{m},\overline{x}')$ we have:
$$\sin d(\overline{m},\overline{y})\sin 
d(\overline{m},\overline{x}')\cos
\overline{\alpha}' = \cos
d(\overline{y},\overline{x}') - \cos
d(\overline{m},\overline{y}).\cos
d(\overline{m},\overline{x}')$$
   $$=\cos (\frac{\pi}{2}-\delta ) - \cos 
(\frac{\pi}{2}-\delta-\epsilon).\cos
d(\overline{m},\overline{x}')$$
$$= \sin (\delta ) -\sin (\delta +\epsilon)\cos 
d(\overline{m},\overline{x}')\le 0$$
which means that  $$\cos 
d(\overline{m},\overline{x}')\ge\frac{\sin(\delta)}{\sin(\delta
+\epsilon)}.$$

It follows that $$L(\sigma) =
d(\overline{m},\overline{x}')\le
\cos^{-1} \left(\frac{\sin(\delta)}{\sin(\delta +
\epsilon)}\right)=\tau(\epsilon)$$ with
$\tau(\epsilon)\mapsto 0$
 when $\epsilon\mapsto 0.$

 Hence $$d(m,x)= d(\overline{m},\overline{x})\ge \pi
-\tau(\epsilon).$$
\subsection{Lemma}\label{111111}

 \textit{Let  $M$ be a complete connected
$n$-dimensional Riemannian 
manifold
with sectional  curvature $ K_{M}\ge 1$ and
$diam(M)\ge \pi-\epsilon,$
 then for all $x\in M$
there exists
a point $x'$ such that
 $d(x,x')\ge \pi-\Gamma(\epsilon)$ hence $Rad(M)\ge
 \pi-\Gamma(\epsilon)$ with $\Gamma(\epsilon)\mapsto
0$
  when $\epsilon\mapsto 0$ where $$Rad(M)= min_{x\in
M}max_{y\in 
M}d(x,y)$$  }

 \begin{center}Proof\end{center}
Since $M$ is compact, its injectivity radius $i(M)$ is positive and
set $r_{0}$  a positive number not larger than $i(M).$ 
  Let $p, q$ be two  points in $M$ such that
$d(p,q) = diam(M)\ge \pi-\epsilon.$ Let $x\in M$ and
suppose
  $x\ne p$ and $x\ne q.$  Consider the triangle $
(p,q,x)$ and let $y$ be a point of the segment $[q,x]$
  such that $d(q,y) = \frac{r}{2}$ with $r < r_{0}$
and
  $d(p,x) >2r$  $d(q,x)>2r.$ For any point
  $s\in B(y,r)$
  the function  $$z\mapsto e_{yz}(s) = d(y,s) + d(z,s)
-d(y,z)$$
is continuous and if $z$
   is on the prolongation of the  geodesic
joining $y$ to $s$ we have: $e_{yz}(s) = 0$ ( this is 
possible
since $r < r_{0}$).

Hence the function $z\mapsto e_{yz}(q)$ is continuous
on the sphere
 $S(y,r)$, and consequently there  exists
   $z\in S(y,r)$ such that
 \begin{equation}\label{101010}
e_{yz}(q) = d(y,q) + d(z,q) -d(y,z)<
\epsilon\end{equation}

For any  point $v\in M$ we have:
$$d(p,v) + d(q,v) + d(p,q)\le 2\pi $$ hence
 \begin{equation}\label{22222}\vert d(p,v) + 
d(q,v)-\pi\vert\le\epsilon.\end{equation}
Let $\overline{p}'$ be the antipodal point of
$\overline{p}$ on the
 sphere $S^{2}.$
$$d(\overline{p},\overline{q}) = d(p,q);$$
$$d(\overline{z},\overline{p}') = \pi-
d(\overline{z},\overline{p} )= 
\pi -d(p,z).$$
We have $$\vert 
d(\overline{q},\overline{y}-d(\overline{y},\overline{p}')\vert
= \vert d(\overline{q},\overline{y})-\pi+ 
d(\overline{p},\overline{y})\vert
< \epsilon.$$
In the same way, we have: $$\vert
d(\overline{q},\overline{z})-\pi+ 
d(\overline{p},\overline{z})\vert = \vert d(q,z)-\pi+
d(p,z)\vert 
<\epsilon$$
Hence\begin{equation} \label{abd}
d(\overline{y},\overline{z}) = d(y,z)
d(q,y) + d(q,z)-\epsilon \ge \pi-\epsilon -d(p,y) 
+d(q,z)-\epsilon\end{equation}
$$ \ge 
\pi-d(\overline{p},\overline{y}) +
d(\overline{q},\overline{z})-2\epsilon
\ge d(\overline{p}',\overline{y}) +
d(\overline{q},\overline{z})-2\epsilon$$
\begin{equation}\label{trois}\ge
d(\overline{p}',\overline{y}) +
d(\overline{p}',\overline{z})-3\epsilon\end{equation}
Thus $d(\overline{p}',\overline{y}) >
\frac{r}{2}-\frac{3}{2}\epsilon$ 
and
$d(\overline{p}',\overline{z}) >
\frac{r}{2}-\frac{3}{2}\epsilon.$

Suppose
$d(\overline{p}',\overline{y}) \le
\frac{r}{2}-\frac{3}{2}
\epsilon$ then
$d(\overline{p},\overline{y}) \ge
\pi-\frac{r}{2}+\frac{3}{2}\epsilon$

which contradicts (\ref{22222}).

Let $$\overline{l}_{1}' =
d(\overline{p}',\overline{y}); 
\overline{l}_{2}' = d(\overline{p}',\overline{z})
\textrm{ and }
\overline{l}_{0}' = d(\overline{z},\overline{y}).$$
Applying
the Topogonov's theorem  to the triangle
$(\overline{y},\overline{p}',\overline{z})$ , we have:
$$\sin\overline{l}_{1}'\sin\overline{l}_{2}'\cos\measuredangle\overline{p}'=\cos\overline{l}_{0}'-\cos\overline{l}_{1}'\cos\overline{l}_{2}'$$
From  inequalty (\ref{trois}) we 
get:$$\sin\overline{l}_{1}'\sin\overline{l}_{2}'\cos\measuredangle\overline{p}'
< 
\cos(\overline{l}_{1}'+\overline{l}_{2}'-
3\epsilon)-\cos\overline{l}'_{1}\cos\overline{l}'{2}$$
$$<-\cos\overline{l}'_{1}\cos\overline{l}'_{2}(1-\cos3\epsilon)-\sin\overline{l}'_{1}\sin\overline{l}'_{2}\cos3\epsilon
-\sin(
\overline{l}'_{1}+ \overline{l}'_{2})\sin 3\epsilon$$
Hence $$ cos\measuredangle\overline{p}' <
-\cos3\epsilon 
+(cotg\overline{l}'_{2}+ cotg\overline{l}'_{1})\sin
3\epsilon
-cotg\overline{l}'_{1}.cotg\overline{l}'_{2}(1-\cos3\epsilon).$$
Thus $\measuredangle\overline{p}'> \pi
-\tau_{1}(\epsilon)$
with $$
\tau_{1}(\epsilon)=\cos^{-1}\left(cos3\epsilon- 
(cotg\overline{l}'_{2}+ cotg\overline{l}'_{1})\sin
3\epsilon
+
cotg\overline{l}'_{1}.cotg\overline{l}'_{2}(1-\cos3\epsilon)\right)$$
and $ \tau_{1}(\epsilon)\mapsto 0$ as
$ \epsilon \mapsto 0.$

The trigonometric  law on the sphere shows that the
angle at
$\overline{p}'$ of the triangle 
$(\overline{y},\overline{p}',\overline{z})$ is equal
to the
angle at
$\overline{p}$ of the triangle 
$(\overline{y},\overline{p},\overline{z})$ which  is
 not bigger than the  angle at $p$ of the triangle
$(y,p,z)$  in $M.$

Hence $\measuredangle (y,p,z) > \pi
-\tau_{1}(\epsilon).$

By applying the relation (\ref{22222}) to  $x$ and $y$
we have:
\begin{equation}\label{22224}-\epsilon\le d(p,x) + 
d(q,x)-\pi\le\epsilon,\end{equation}
\begin{equation}\label{22225}-\epsilon\le d(p,y) + 
d(q,y)-\pi\le\epsilon.\end{equation}

Since $y\in[q,r],$
 we  conclude from (\ref{22224}) and (\ref{22225})
that
$e_{py}(x)< 2\epsilon$, which shows that the angle at
$x$ of the 
triangle $(p,x,y)$ is close to
 $\pi$ and
  consequently its  angle at $p$ is small. there exists

$\tau_{2}(\epsilon)$ such that $\measuredangle p\le
\tau_{2}(\epsilon)$
  where $\tau_{2}(\epsilon)\mapsto 0$ as
$\epsilon\mapsto 0.$

 Take $r$ small enough; then $d(p,x)< d(p,y);$
otherwise $d(p,x)$ is 
close to $\pi$ and we conclude by taking
  $x' = p.$\

 Let  $\tilde{x}\in [p,y]$ such that  $d(p,\tilde{x})
= d(p,x);$ then

  \begin{equation}\label{11221} d(x,\tilde{x})\le
 \pi\measuredangle
(y,p,x)=\tau_{3}(\epsilon).\end{equation}

Since
$$d(p,z) + d(z,q)\ge d(p,q)\ge \pi -\epsilon$$ and
 $$d(y,z)= r > d(q,y) +d(q,z)-\epsilon$$ we have:

if  $d(p,x)\le d(q,z)+\epsilon$ then  $r < \frac{2}{3}\epsilon$ and it suffices to take $x' = q$ and $\Gamma(\epsilon)\le \frac{7}{3}\epsilon ;$
 if $d(p,x) > d(q,z)+\epsilon$ then
$$d(p,z) \ge\pi-d(q,z)-\epsilon > \pi - d(p,x)$$
hence there exists a point $x'\in[p,z]$ such that
$$d(p,x') = \pi- d(p,x).$$

It suffices to show that $d(x',\tilde{x})\ge
\pi-\tau_{4}(\epsilon).$

Applying the  Toponogov theorem to the triangle
$(\tilde{x},p,x'),$
we get

$$\cos d(x',\tilde{x}) = \cos d(p,x')\cos
d(p,\tilde{x}) + \sin 
d(p,x')\sin d(p,\tilde{x})
\cos\measuredangle p$$
$$ = -\cos^{2} d(p,\tilde{x}) + \sin^{2} 
d(p,\tilde{x})\cos\measuredangle p$$
Since $x'\in[p,z]$ and $\tilde{x}\in[p,y]$  the angle
at
$\overline{p}$ of the triangle 
$(\overline{y},\overline{p},\overline{z})$ is less or
equal to the
angle at $\overline{p}$ of the triangle
$(\overline{\tilde{x}},\overline{p},\overline{x'})$
which is not bigger 
than the angle at $p$
of triangle $(\tilde{x},p,x') $ in $M.$

$$\measuredangle(\tilde{x},p,x')\ge 
\measuredangle(\overline{\tilde{x}},\overline{p},\overline{x'}\ge
\measuredangle(\overline{y},\overline{p},\overline{z})\ge
\pi-\tau_{4}(\epsilon).$$
Hence
$$ cos d(x',\tilde{x}) \le -\cos^{2} d(p,\tilde{x}) +
\sin^{2} 
d(p,\tilde{x})\cos(\pi-\tau_{4}(\epsilon)$$
$$=-\cos (\tau_{4}(\epsilon)) - \cos^{2} 
d(p,x)\left(1-\cos(\tau_{4}(\epsilon))\right)\le
-\cos(\tau_{4}(\epsilon))$$
$$\Rightarrow d(\tilde{x},x') > \pi-
\tau_{4}(\epsilon).$$

 From the triangle inequality and the inequality
(\ref{11221}) we have
  $$d(x,x')\ge d(\tilde{x},x')-d(\tilde{x},x)\ge
\pi-\tau_{5}$$
Thus,  lemma \ref{111111} follows.

In [Gr] M. Gromov generalized the classic notion of
Hausdorff
distance  between two compact subsets of the same
metric space.
He considered the set of compact Riemannian manifolds
as
 a subset of the set of all compact metric
spaces.

\subsection{Definitions}

\textit{1)\quad Let $ X,Y$ be two metric espaces; a
map $f: X 
\longrightarrow Y$ is said to be an
 $\epsilon$-approximation if the image set $f(X)$ is
$\epsilon$-dense 
in $Y$ and, for  any $x,y \in X,$
 $\vert d(f(x),f(y))-d(x,y)\vert < \epsilon.$}

\textit{2)\quad The  Gromov-Hausdorff distance
   $d_{GH}(X,Y)$ between  $X$ and  $Y$ is the infimum
of values of 
$\epsilon > 0$
 such that there exist  $\epsilon$- approximations 
$f: X 
\longrightarrow Y$
 and $g: Y \longrightarrow X.$}
In [Co1] and [Co2] Colding showed the two equivalent conditions:

1) $Rad(M)\ge \pi-\epsilon$ 

2)
$d_{GH}\left(M , 
S^{n}(1)
\right)\le \tau_{5}(\epsilon)$ with
$\tau_{5}(\epsilon)\mapsto 0$
as $\epsilon\mapsto 0.$

Then theorem \ref{**} follows from these conditions and the theorem A.

\section{references}
[C] J.Cheeger, Critical points of distance functions
and applications 
to
geometry , Lectures notes 1504 (1991) 1-38.
\newline
[CC] J.Cheeger and T. Colding, On the structure of
spaces with Ricci
curvature bounded from below, J. Diff. Geom. 46 (1997)
406-480.
        \newline
[Ch] S.Y.Cheng, Eigenvalue Comparison Theorem and
Geometric 
Applications
 Math.Z.143(1975) p. 289-297.
\newline
[Co1] T.H. Colding, Shape of manifolds with positive Ricci curvature, Invent. Math. 124 (1996), 175-191
\newline
[Co2] T.H. Colding, Large manifolds with positive Ricci curvature, Invent. Math. 124 (1996), 193-214
\newline
 [Gr] M.Gromov, Structures m\'etriques  pour  les 
vari\'et\'es  
Riemanniennes,
 R\'edig\'e par J.Lafontaine et P.Pansu,cedic 1981.
\newline
[GS] K. Grove and K. Shiohama, A generalized sphere
theorem, Ann. Math.
106 (1977), 201-211.
\newline
[Pe] P. Pertersen, Comparison geometry problem list,
Riemannian 
geometry (Waterloo, ON,(1993) 87-115,
fields Inst. Monogr., 4, Amer. Math. Soc., Providence,
RI 1996.
        \newline
[S] K. Shiohama, A sphere theorem for manifolds of
positive Ricci 
curvature, Trans. Amer. Math. Soc.
275 n.2 (1983) 811-819.
\newline
[Sh] Z. Shen , Complete manifolds with nonnegative
Ricci curvature and
large volume growth, Invent. Math. 125 (1996) 393-404.
                \newline
[SS] J. Sha and Z. Shen, Complete manifolds with
nonnegative curvature 
and quadratically nonnegative curved infinity, Amer.
J. Math. 119 (1997) 1399-1404.
                      \newline
[T] V.A.Toponogov, Computation of length of a closed
geodesic on convex 
surface, Dokl. akad.
 Nank SSSR 124 (1959) 282-284.
\newline 
[Xi1] C.Y. Xia , Large volume growth and topology of
open manifolds,
 Math. Z. 239 (2002) 515-526.
\newline
[Xi2] "  " , Complete manifolds with sectional
curvature bounded below
and large volume growth, Bull. London Math. Soc. 34
(2002) 229-235.
                   \newline
[Xi3] " " , Some applications of critical point theory
of distance 
functions on Riemannian manifolds, 
Compos.math.132, n.1(2002)pp 49-55. 

\vspace{1cm}
Bazanfar\'e Mahaman\\
Universit\'e Abdou Moumouni \\
D\'epartement de Math\'ematiques et Informatique Niamey- Niger\\
E-mail: bmahaman @yahoo.fr

    \end{document}